\documentclass[reqno]{amsart}
\usepackage{graphicx}
\usepackage{amsmath,amssymb}

\newtheorem{thm}{Theorem}[section]
\newtheorem{cor}[thm]{Corollary}
\newtheorem{lem}[thm]{Lemma}
\newtheorem{prop}[thm]{Proposition}
\newtheorem{rem}[thm]{Remark}

\newcommand{\To}{\longrightarrow}

\begin{document}
\large 
\noindent  
  \textbf{Order of approximation in the central limit theorem for associated random variables and a moderate deviation result}\\[0.6 cm]
  M. Sreehari \\
 \textit{ 6-B, Vrundavan Park,
  New Sama Road,
  Chani Road P.O.\\
  Vadodara, 390024, India.}\\[0.5 cm] 
\begin{footnotesize}
\makeatletter{\renewcommand*{\@makefnmark}{}
	\footnotetext{
		E-\textit{mail addresses}: msreehari03@yahoo.co.uk\\	orcid.org/0000-0002-2381-2891}}

\end {footnotesize}

 \textbf{Abstract.}\\
 An estimate of the order of approximation in the central limit theorem for strictly stationary associated random variables with finite moments of order $q > 2$ is obtained.  A moderate deviation result is also obtained. We have a refinement  of recent results in \c{C}a\v{g}in \textit{et al.} (2016).   The order of approximation obtained here is an improvement over the corresponding result in Wood (1983). 
 \section{Introduction}
 A set of random variables (rvs) $\left\{X_1, X_2, \ldots, X_k\right\}$ is said to be associated if for each pair of coordinatewise nondecreasing functions $f, g: R^k\rightarrow R$	
 $$	Cov(f(X_1, X_2, \ldots, X_k), g(X_1, X_2, \ldots, X_k))\geq 0 $$ whenever the covariance exists.\\	A sequence $\left\{X_n\right\}$ of rvs is associated if for every $n \in N$ the family $X_1, X_2, \ldots, X_n$ is associated.\\\\
 In this paper we consider   a strictly stationary sequence of centered square integrable associated rvs $\left\{X_n\right\}$. Central limit theorem (CLT) for  $\left\{X_n\right\}$ was proved by Newman (1980) and a Berry-Ess\'{e}en type theorem giving an estimate of the order of approximation in the  CLT was proved by Wood (1983). In the case of finite third absolute moment $E|X_1|^3$ Wood's result gives an estimate   of the order $O(n^{-1/5})$.  Birkel (1988) obtained a rate of the order $O(n^{-1/2}\log ^2 n)$ under the strong additional assumption  that the Cox-Grimmett coeficients $u(n)$ decrease exponentially. Birkel also provided an interesting example to show the reasonableness of the assumptions to obtain the above order of approximation.  In that example he showed that the above rate cannot be obtained if $u(n)$ decreases  only as a power. Thus there is a huge gap between the results of Wood and Birkel.
  In a recent paper  \c{C}a\v{g}in \textit{et al.} (2016) obtained another estimate of the order of approximation  in the CLT for associated rvs and also obtained a moderate deviation type result. However their estimate in the case of finite third absolute moment $E|X_1|^3$ is quite complicated. \\\\
	We recall here that Wood (1983) mentioned some examples of stationary associated random variables in models for ferromagnets in mathematical physics and in particular Ising model. Further,  large deviation probability and moderate deviation probability investigations received much attention due to their importance in statistical inference and applied probability. We refer to monographs by Vardhan (1984), Dembo and Zeitouni (1998) and Hollander (2000) and recent papers by Wang (2015) and  \c{C}a\v{g}in \textit{et al.} (2016) for other references. These investigations are also useful in the construction of certain counter examples (see, for example, Tikhomirov, 1979 and Birkel, 1988).\\\\  
 We give an  estimate of the order of approximation in the CLT which is a refined version of the result in \c{C}a\v{g}in \textit{et al.} (2016) and also prove corresponding moderate deviation result.  In the case of $X_n$ with finite third absolute moment, when Cox-Grimmett coefficients $ u(n)$ are of order $n^{-\delta}$, the order of approximation in the CLT is proved to go to zero as $n^{-3/8}$ as $\delta \rightarrow \infty$. The main steps in the proof  are the classical decomposition of the partial sum $S_n=\sum_{j=1}^n X_j$ into blocks ( of size $p_n=[n^{1-\alpha}], 0 < \alpha < 1$ ) , coupling them with blocks variables with the same distributions but independent and use the inequality due to Newman (1980). Our approach is similar to that in \c{C}a\v{g}in \textit{et al.} However the estimate of the order of approximation we obtain does not depend on the value of $\alpha$ whereas the same obtained by  \c{C}a\v{g}in \textit{et al.} depends on $\alpha$.  The refinement is in terms of the assumptions, bound and simplification of the steps. This helps us  to get moderate deviation type result too under assumptions milder than those in \c{C}a\v{g}in \textit{et al.}  (2016) and also get an order of approximation in the CLT  which is an improvement over the corresponding result in Wood (1983).\\\\
 The paper is organized as follows. In Section 2 we introduce notation and give some lemmas. In Section 3 we shall have  a set of propositions that will be used  in later sections. Order of approximation in the CLT is investigated in Section 4. Finally a moderate deviation type result is discussed in Section 5.
 \section {Notation}
 Let $\{X_n\}$ be  a strictly stationary sequence of centered square integrable associated rvs.  Set  $E(X_1^2)=\sigma_1^2,\; c_j= Cov(X_1, X_{1+j}),\; S_n=\sum_{j=1}^nX_j,\; ES_n^2=s_n^2$ and $\sigma^2=\sigma_1^2+2\sum_{j=1}^\infty c_j > 0.$\\	
 We assume that  $\sum_{j=1}^\infty c_j < \infty.$ Then	
 $$
 \frac{S_n}{\sigma \sqrt{n}}\stackrel{D}{\rightarrow} Z_1 \sim N(0, 1) 
 $$
 where $N(0, 1)$ denotes the standard normal distribution. The standard proof of this result  involves writing $S_n$ as the sum of blocks of fixed size, approximating the distribution of $S_n$ by the distribution of corresponding sum of coupling block rvs (to be defined shortly) and appealing to the CLT for the coupling block rvs. We need more notation to explain this.
 Define initial blocks $$ Y_{j, n}=\sum_{i=(j-1)p_n+1}^{jp_n} X_i,\;\;\; j=1, 2, \dots, m_n$$
 and $$Y_{m_n+1, n}=\sum_{i=m_np_n+1}^n X_i$$
 where $m_n=[n/p_n],\; p_n < n/2 $ and $[r]$ denotes the largest integer $ \leq r$. Clearly $$ S_n=\sum_{j=1}^{m_n}Y_{j, n}+Y_{m_n+1, n}.$$  
We note that $Y_{j, n}, j=1, 2, \ldots, m_n$ are identically distributed. Further $n- m_n p_n \leq p_n$. We next define  independent coupling blocks $Y_{j, n}^*, j=1, 2, \cdots, m_n$, where $Y_{j, n}^* \stackrel{D}{=} Y_{j, n}.$  Note that since the $X_k$ are strictly stationary, the rvs $Y_{j, n}^*$ are independent and identically distributed.\\Set $p_n=[n^{1-\alpha}]$ where $0< \alpha < 1$.\\
 In what follows limits are taken as $n \rightarrow \infty$ and statements hold for sufficiently large values of $n$.	We make some  of the following assumptions on the covariances $c_j$ and moments of $X_k$ in the remainder:\\
 Assumption $A_1$:  $E|X_k|^q < \infty $ for some  $q > 2 $.\\
 Assumption $A_2$:  $|\frac{s_n^2}{n \sigma^2}-1|=O(n^{-\theta})$ for some $\theta > 0$ where $ s_n^2=ES_n^2$. \\
 Assumption $A_3$:  $u(n)=\sum_{j=n}^\infty c_j < C_1 n^{-\delta}$, where $ \delta > 0$.
 \begin {rem} \label {R1}
 (i)	If $ k_n \rightarrow \infty$ such that $ \frac{k_n}{n}\rightarrow 0$ then the assumption $A_2$ implies  $|\frac{s_n^2}{n \sigma^2} - \frac{s_{k_n}^2}{k_n \sigma^2}| = O( k_n^{-\theta}).$\\
 (ii) By the assumption that $\sum_{i=1}^\infty c_i < \infty$ and the fact that $\sigma^2 - \frac{s_n^2}{n} =2u(n)+\frac{2}{n}\sum_{i=1}^{n-1} j c_j$ it follows that if the assumption $A_2$ holds for some $\theta > 0$ then the assumption $A_3$ holds for $\delta = \theta$ and conversely.\\
 (iii) Under the  assumption $A_1$ there exist positive constants $A$ and $B$ such that for all the positive integers $n$, $A \;n^{1/2} < s_n < B\; n^{1/2}$ and $A \;n^{q/2} < E|S_n|^q < B\; n^{q/2}.$ (see (2.16) in Birkel, 1988)
\end{rem}
Here and elsewhere $C_1, C_2, \dots$ are positive constants independent of $n$. Further  $\eta_1, \eta_2, ...$ are constants with absolute values $ \leq 1$. The following result is known.
\begin {lem}(Newman's inequality, 1980)  \label {L1}
Suppose $U_1, U_2, \cdots, U_n$ are associated rvs with finite variances. Then for any real numbers $t_1, t_2, \cdots, t_n$
$$\left|E\left(\exp^{i\sum_{j=1}^nt_j\;U_j}\right)-\prod_{j=1}^n E\left(\exp^{it_j\;U_j}\right)\right| \leq \sum_{i=1, j>i}^n|t_i||t_j|Cov(U_i, U_j).$$
\end {lem}
\begin{rem} \label {R0}
If $\{X_n\}$ is a sequence of associated rvs, the block rvs $Y_{1, n}, Y_{2, n}, \cdots,$ $ Y_{m_n, n}$ are associated. Further the characteristic functions satisfy
\begin {equation} \label{E0}
E\left( e^{i\sum_{j=1}^{m_n}t_jY_{j, n}^*}\right)=\prod_{j=1}^{m_n}E\left( e^{it_jY_{j, n}^*}\right)= \prod_{j=1}^{m_n}E\left( e^{it_jY_{j, n}}\right)
\end{equation}
because $Y_{j, n}^*$ are independent and $Y_{j, n}^*\stackrel{D}{=}Y_{j,n}$.
\end{rem}
Let $T_1=a_n n^{\alpha/2} $ and $T_2= b_n n^{\alpha/2}$ where $a_n =(\log n)^a$ and $b_n=(\log n)^b$ with $a < b < 0.$ 
\section{Some preliminary results}
In this section we discuss some preliminary results that will be used   later and   these are of independent interest too.	
The following result notes that while dealing with asymptotic properties of $S_n/(\sigma s_n)$ it is adequate to consider the sum  $\sum_{j=1}^{m_n} Y_{j, n}$. 
\begin{prop}  \label {P1}
Suppose the assumption $A_1$  holds. Then	for $\mu_n=n^{-3\alpha/8},\;\;\; 0 < \alpha < 1$ 
$$P(|Y_{m_n+1, n}| > \mu_n s_n) < C_2n^{-q\alpha/8}.$$
\end{prop}
To see this note that because of stationarity of $\{X_n\}$ $$P(|Y_{m_n+1, n}| > \mu_n s_n) <\frac{E|S_{n-m_np_n}|^q}{\mu_n^q s_n^q}.$$
The result follows now from Remark \ref{R1} and the assumption $A_1$. 

\begin {rem} This is an improvement of the result in Step 3  of the Theorem 3.1 in \c{C}a\v{g}in \textit{et al.} (2016) as it does not put any restriction on $\alpha$ and $q$.\\(ii) One can chose $\mu_n=n^{-\mu\alpha}, \;\;0 < \mu <1/2$  but the calculations become too complicated. See Remark 4.2 below.
\end {rem}	
Next we approximate the distribution of the sum of the original rvs by that of the coupling blocks; i.e., the distribution of $\sum_{j=1}^{m_n} Y_{j, n}$ by that of $\sum_{j=1}^{m_n} Y_{j, n}^*.$ The method of approximation is based on the celebrated Berry-Ess\'{e}en inequality and Newman's inequality for associated rvs. 

\begin{prop} \label {P2}. Suppose the assumptions $A_1$ and $A_2$ hold. Then 
$$\sup_{x \in R} \left|P\left(\sum_{j=1}^{m_n} Y_{j, n} \leq x s_n\right) - P\left(\sum_{j=1}^{m_n} Y_{j, n}^* \leq x s_n\right)\right|$$ $$ <  C_3\frac{b_n^2}{n^{\theta-\alpha(1+\theta)}}I\left(\frac{2\theta}{3+2\theta} \leq \alpha < \frac{\theta}{1+\theta}\right)+C_4\frac{1}{b_n n^{\alpha/2}}I\left(\alpha < \frac{2\theta}{3+2\theta}\right).$$
\end{prop}	
\textbf{Proof}
By the Berry- Ess\'{e}en inequality and (\ref{E0}) we have
$$\sup_{x \in R} \left|P\left(\sum_{j=1}^{m_n} Y_{j, n} \leq x s_n\right) - P\left(\sum_{j=1}^{m_n} Y_{j, n}^* \leq x s_n\right)\right|$$
$$ < C_5\int_{-T_2}^{T_2}\frac{1}{|t|}\left|E\left(e^{i\frac{t}{s_n}\sum_{j=1}^{m_n} Y_{j, n}}\right)-\prod_{j=1}^{m_n}E\left(e^{i\frac{t}{s_n}Y_{j, n}^*}\right)\right|dt + \frac{C_6}{T_2}$$

\begin {equation}\label {E1}
=  C_5\int_{-T_2}^{T_2}\frac{1}{|t|}\left|E\left(e^{i\frac{t}{s_n}\sum_{j=1}^{m_n} Y_{j, n}}\right)-\prod_{j=1}^{m_n}E\left(e^{i\frac{t}{s_n}Y_{j, n}}\right)\right|dt + \frac{C_6}{T_2}.
\end{equation}
By the Lemma \ref {L1} with $ U_j=Y_{j, n}, \; j=1, 2,\cdots, m_n$ we have  
$$	\left|E\left(e^{i\frac{t}{s_n}\sum_{j=1}^{m_n} Y_{j, n}}\right)-\prod_{j=1}^{m_n}E\left(e^{i\frac{t}{s_n}Y_{j, n}}\right)\right|
\leq \frac{t^2}{s_n^2}\sum_{j, k=1, j > k}^{m_n} Cov(Y_{j, n}, Y_{k, n}) 
$$
$$
= \frac{t^2 m_np_n}{2 s_n^2}\left|\frac{s_{m_np_n}^2}{m_n p_n} - \frac{s_{p_n}^2} {p_n}\right|. 
$$
In view of the Remark \ref {R1}  we then have from  (\ref{E1})
$$	\sup_{x \in R} \left|P\left(\sum_{j=1}^{m_n} Y_{j, n} \leq x s_n\right) - P\left(\sum_{j=1}^{m_n} Y_{j, n}^* \leq x s_n\right)\right|<\frac{C_7}{p_n^\theta}  \int_{-T_2}^{T_2}|t|dt +\frac{C_6}{T_2}$$
$$= \frac{C_7\;T_n^2}{p_n^\theta}+\frac{C_6}{T_2}.$$
Recalling that $p_n=[n^{1-\alpha}]$ , $T_2= b_n\; n^{\alpha/2}$ we note that the right side above goes to zero only for $\alpha <\theta/(1+\theta)$. Further $(1-\alpha)\theta \leq 3\alpha/2$ if and only if $\alpha  \geq 2\theta/(3+2\theta).$ Hence 
$$\sup_{x \in R} \left|P\left(\sum_{j=1}^{m_n} Y_{j, n} \leq x s_n\right) - P\left(\sum_{j=1}^{m_n} Y_{j, n}^* \leq x s_n\right)\right|$$ $$ \leq C_6\;\frac{1}{b_n n^{\alpha/2}}\;I\left(\alpha < \frac{2\theta}{3+2\theta}\right) + C_8\; \frac{b_n^2}{n^{(1-\alpha)\theta-\alpha}}\;I\left(\frac{2\theta}{3+2\theta} \leq \alpha < \frac{\theta}{1+\theta}\right).$$
This completes the proof of the Proposition \ref{P2}.
\begin{rem}
In \c{C}a\v{g}in \textit{et al.} (2016) the above bound was obtained separately for the odd numbered blocks and the even numbered blocks.  Further, the bound obtained above goes to zero faster than their corresponding bound.  
\end{rem}
Our next result is concerned with the approximation of the characteristic function of the sum of coupling blocks by the characteristic function of an appropriate  normal variable.
\begin{prop} Denote $\varphi _j(t)= E\left(e^{itY_{j, n}^*}\right)$. Then under the assumptions $A_1$ and $A_2$, for $|t| < T_2$ 
$$\left|\prod_{j=1}^{m_n}\varphi_j(\frac{t}{s_n})-e^{-\frac{m_nt^2 s_{p_n}^2}{2s_n^2}}\right| \leq  C_9\;\frac{m_n|t|^q  p_n^{q/2}}{s_n^q}e^{-\frac{m_nt^2 s_{p_n}^2}{2s_n^2}}.$$ \label{P3}	
\end{prop}	
\textbf{Proof} Let us first consider the case  $2 < q < 3$.
Note that since $Y_{j, n}^*\stackrel{D}{=}Y_{j, n}$ $$\varphi_j(t/s_n)=1-\frac{t^2s_{p_n}^2}{2s_n^2}+\eta_1\;\frac{|t|^q}{q(q-1)s_n^q} E|Y_{j, n}|^q.	$$ 
For $|t|< T_1= a_n \;n^{\alpha/2}$, with $ a_n=(\log n)^a,\; a < 0$
$$\frac{t^2\;s_{p_n}^2}{s_n^2} < C_{10} \;a_n^2 \rightarrow 0.$$
Further	
$$ \frac{|t|^q}{s_n^q} E|Y_{j, n}|^q < C_{11}\; a_n^q \rightarrow 0.$$
Hence $	|\varphi_j(t/s_n)- 1| \rightarrow 0$ and therefore 	$\varphi_j(t/s_n)$ is bounded away from 0 for $|t| < T_1$  so that we can take its logarithm. Then for each $j$
$$ \log \varphi_j(t/s_n) = -\frac{t^2s_{p_n}^2}{2s_n^2}+\eta_1\;\frac{|t|^q}{q(q-1)s_n^q} E|Y_{j, n}|^q + \eta_2\;\left[-\frac{t^2s_{p_n}^2}{2s_n^2}+\eta_1\;\frac{|t|^q}{q(q-1)s_n^q} E|Y_{j, n}|^q\right]^2$$
$$ = -\frac{t^2s_{p_n}^2}{2s_n^2}+\eta_3\;\frac{|t|^q}{q(q-1)s_n^q} E|Y_{j, n}|^q.$$
Then using the fact $|e^x-1|< |x|\; e^{|x|}$ we get  $$\left|\prod_{j=1}^{m_n}\varphi_j(t/s_n)- e^{-\frac{m_n\;t^2 s_{p_n}^2}{2s_n^2}}\right| \leq e^{-\frac{m_nt^2 s_{p_n}^2}{2s_n^2}}\;\frac{m_n|t|^qE|Y_{j, n}|^q}{2\;s_n^2}\;e^{\frac{m_n|t|^qE|Y_{j, n}|^q}{s_n^q}}$$
Note that $\frac{|t|^{q-2}E|Y_{1, n}|^{q}}{s_n^{q-2}\;s_{p_n}^2} < a_n^{q-2} \rightarrow 0$ so that $\frac{m_n|t|^qE|Y_{1, n}|}{s_n^q} < \frac{m_n t^2s_{p_n}^2}{4\;s_n^2} $ and hence
\begin{equation}
\left|\prod_{j=1}^{m_n}\varphi_j(t/s_n)- e^{-\frac{m_n\;t^2 s_{p_n}^2}{2s_n^2}}\right| < C_{12}\; \frac{m_n\;|t|^q\; E|Y_{1, n}|^q }{s_n^q} \;\;e^{-\frac{m_n\;t^2 \;s_{p_n}^2}{4\;s_n^2}}\label{E2}
\end{equation}
for $|t| < T_1$.  We shall prove that the relation (\ref {E2}) holds for $ T_1 \leq |t| < T_2$ also.\\\\
Let $W_j,\; j=1, 2, \cdots, m_n$ be rvs such that for each $j,\; W_j$ is independent of $Y_{j, n}^*$ and distributed as $Y_{j, n}^*$. Then $E(W_j-Y_{j, n}^*)=0, E(W_j-Y_{j, n}^*)^2 =2s_{p_n}^2$ and $E|W_j-Y_{j, n}^*|^q \leq 2^q E|Y_{j, n}^*|^q$. Further $$|\varphi_j(t/s_n)|^2= E\left(e^{i\frac{t}{s_n}(W_j-Y_{j, n}^*)}\right)$$
$$=1-\frac{t^2s_{p_n}^2}{s_n^2}+ \theta_4\frac{2^q |t|^qE|Y_{j, n}^*|^q}{q(q-1)s_n^q}.$$
Note that for $|t|< T_2= b_n\; n^{\alpha/2}$ by the Lemma \ref{L1}, $$\left|\theta_4\frac{2^q |t|^qE|Y_{j, n}^*|^q}{q(q-1)s_n^q}\right| < C_{13}\frac{t^2 s_{p_n}^2}{s_n^2}\left[\frac{T_2}{n^{\alpha/2}}\right]^{q-2}$$
$$ < C_{14}\frac{t^2 s_{p_n}^2 b_n^{q-2}}{s_n^2} < \frac{3t^2 s_{p_n}^2}{4s_n^2}$$
since $b_n\rightarrow 0$. Hence $$|\varphi_j(t/s_n)|^2 < 1-\frac{t^2s_{p_n}^2}{4s_n^2}.$$
Since $\frac{t^2s_{p_n}^2}{s_n^2}\rightarrow 0$, using the fact $1-u  < e^{-u}$ for $u > 0$  we have
$|\varphi_j(t/s_n)|^2 < \exp\left(-\frac{t^2s_{p_n}^2}{4s_n^2}\right).$ 
Thus for $|t| < T_2$
\begin {equation}
\left|\prod_{j=1}^{m_n}\varphi_j(t/s_n)- e^{-\frac{m_n\;t^2s_{p_n}^2}{2s_n^2}}\right|<2 e^{-\frac{m_n\;t^2s_{p_n}^2}{4\;s_n^2}}.\label {E3}
\end{equation}
Now to complete the proof of the claim  that (\ref{E2}) holds for $T_1 \leq |t| < T_2$ also, consider 
$$C_{12}\frac{m_n|t|^qE|Y_{1, n}^*|^q}{s_n^q} > C_{12}\frac{m_nT_1^q [ES_{p_n}^2]^{q/2}}{s_n^q} > C_{15}n^\alpha a_n^q \rightarrow \infty.$$
Hence for $n$ large  $$ C_{12}\frac{m_n|t|^qE|Y_{1, n}^*|^q}{s_n^q} > 2, $$
and the claim that (\ref{E2}) holds for $T_1 < |t| < T_2 $ also follows from (\ref {E3}) \\

The result of the Proposition then follows from (\ref{E2}) and the Remark  \ref{R1} in the case $2 < q < 3$.\\
In the case $q \geq 3$ we can expand $\log \varphi_j(t/s_n)$  using the third moment also and similar calculations lead to the same bound as above and hence the  Proposition holds true for $q \geq 3.$
\begin{rem} The above proof  is similar to that  in the Theorem 4.1 in \c{C}a\v{g}in \textit{et al.} (2016) but has greater clarity. Further the  final bound is a bit different because we use different values of $T$s. 
\end{rem}
\begin {cor} \label {C1}
Suppose the assumptions $A_1$ and  $A_2$  hold. Then 
$$\int_{-T_2}^{T_2}\frac{1}{|t|}\left|\prod_{j=1}^{m_n}\varphi_j(\frac{t}{s_n})-e^{-\frac{m_nt^2 s_{p_n}^2}{2s_n^2}}\right|dt \leq C_{16}\; \frac{1}{n^{\alpha(q-2)/2}}$$
in the case $2 < q < 3$ and the above inequality holds with $q=3$ giving the bound $C_{16}n^{-\alpha/2}$ in the case $q \geq 3.$
\end{cor}
Here we use the fact that the normal distribution has finite moments.\\
The final result of this section is to approximate the normal distribution with the characteristic function  $ e^{-\frac{m_nt^2 s_{p_n}^2}{2s_n^2}}$ by the standard normal distribution.
\begin {rem} The proof here is essentially the  same as that of the  Theorem 4.1 in \c{C}a\v{g}in \textit{et al.} (2016) but is included for completeness.
\end{rem}
\begin {prop} \label {P4}
Let $G_n(x)$ be the distribution function with the characteristic function $ exp(-\frac{m_nt^2 s_{p_n}^2}{2s_n^2})$ and $\Phi$ be the standard normal distribution function. Then $$\sup_{x \in R} |G_n(x) - \Phi (x)| < C_{17}\; \frac{1}{b_n\;n^{\alpha/2}}\;I\left(\alpha \leq \frac{2\theta}{1+2\theta}\right)$$
$$+ C_{18}\; \frac{1}{n^{(1-\alpha)\theta}}\;I\left(\alpha > \frac{2\theta}{1+2\theta}\right).$$
\end{prop} 
\textbf{Proof} By the Berry-Ess\'{e}en inequality
$$\sup_{x \in R}|G(x)-\Phi(x)| \leq C_{19} \;\int_{-T_2}^{T_2}\frac{1}{|t|}\left|e^{-\frac{m_nt^2 s_{p_n}^2}{2s_n^2}}-e^{-t^2/2}\right| dt + C_{20}\;\frac{1}{T_2}.$$
Using again the fact that $|e^a-1| \leq |a| e^{|a|} $ and recalling that $\frac{m_n\;s_{p_n}^2}{s_n^2} \rightarrow 1$ we have for large $n$  $$\frac{1}{|t|}\left|e^{-\frac{m_nt^2 s_{p_n}^2}{2s_n^2}}-e^{-t^2/2}\right| \leq e^{-\frac{t^2}{2}}\left|\frac{t}{2}\left(\frac{m_n s_{p_n}^2}{s_n^2}- 1\right) \right|e^{\frac{t^2}{2}\left|\frac{m_n s_{p_n}^2}{s_n^2}- 1\right|}$$
$$\leq \left|\frac{m_n s_{p_n}^2}{s_n^2}- 1\right| \frac{|t|}{2}e^{-\frac{t^2}{4}}.  $$
Since the normal distribution has all the moments finite 
$$\sup_{x \in R}|G(x)-\Phi(x)| \leq C_{19}\left|\frac{m_n s_{p_n}^2}{s_n^2}- 1\right|+C_{20}\frac{1}{T_2}.$$
$$ < C_{21} \frac{1}{n^{(1 -\alpha)\theta}}+C_{20}\frac{1}{n^{\alpha/2}b_n}.$$	
Note that $n^{-(1-\alpha)\theta} \rightarrow 0$  faster than $ b_n^{-1}\; n^{-\alpha/2}$ for $\alpha < \frac{2\theta}{1+2\theta}$ while for  $\alpha > \frac{2\theta}{1+2\theta}$   $ b_n^{-1}\; n^{-\alpha/2}\rightarrow 0$ faster than $n^{-(1-\alpha)\theta}$ giving us the stated bound.
\begin{rem}The bound obtained here is a better bound than the corresponding bound in \c{C}a\v{g}in \textit{et al.} (2016).
\end{rem}
\section{Order of approximation in the CLT }
We now obtain an  estimate of the order of approximation in the CLT which is a refined version of the result in \c{C}a\v{g}in \textit{et al.} (2016).  The refinement is in terms of the assumptions, bound and simplification of the steps. It also provides a better bound than the bound of order $n^{-1/5}$ obtained from  Wood's result under the assumption of finiteness of third absolute moments. See Corollary 4.14 in Oliveira (2012).
\begin{thm} \label{T1}
Lett the assumptions $A_1$ and $A_2$ hold. Then
$$\sup_{x \in R} |P(S_n \leq x s_n) - \Phi (x)| < C_{22} [n^{-\frac{\theta(q-2)}{q+2\theta}}I(2 < q \leq 8/3 )$$
$$ + n^{-\frac{q\theta}{q+8+8\theta}}I(8/3 \leq q < 3)+ \;n^{-\frac{3\theta}{11+8\theta}}I(q \geq 3)].$$
In particular when $q=3$ the bound becomes $C_{22} n^{-\frac{3\theta}{11+8\theta}}$.
\end{thm}
\textbf{Proof}
Recall  $\mu_n= n^{-3\alpha/8}. $ Then by the Proposition \ref{P1}, after making  elementary adjustments, we get
$$\sup_{x \in R} |P(S_n \leq x s_n) - \Phi (x)| \leq  \sup_{x \in R}\left|P\left (\sum_{j=1}^{m_n}Y_{j, n} \leq x s_n\right)- \Phi (x) \right|$$
$$+P(|Y_{m_n+1, n}| > \mu_n s_n)
 + \sup_{x \in R} |\Phi (x+\mu_n)- \Phi(x)|$$
$$
< \sup_{x \in R}\left|P\left (\sum_{j=1}^{m_n}Y_{j, n} \leq x s_n\right)- \Phi (x)\right|+ C_{23}\frac{1}{n^{q\alpha/8}} + C_{24}\frac{1}{n^{3\alpha/8}}$$
\begin {equation}\label{E4} 
< \sup_{x \in R}\left|P\left (\sum_{j=1}^{m_n}Y_{j, n} \leq x s_n\right)- \Phi (x)\right|+ C_{23}\frac{1}{n^{q\alpha/8}}+ C_{24}\frac{1}{n^{3\alpha/8}}. 
\end{equation}
 Further by the Berry - Ess\'{e}en inequality
$$\sup_{x \in R}\left|P\left (\sum_{j=1}^{m_n}Y_{j, n} \leq x s_n\right)- \Phi (x)\right|
\leq \sup_{x \in R}\left|P\left (\sum_{j=1}^{m_n}Y_{j, n} \leq x s_n\right)- P\left (\sum_{j=1}^{m_n}Y_{j, n}^* \leq x s_n\right)\right|$$
$$ + \sup_{x \in R}\left|P\left (\sum_{j=1}^{m_n}Y_{j, n}^* \leq x s_n\right)-  \Phi (x)\right|$$
$$ <\sup_{x \in R}\left|P\left (\sum_{j=1}^{m_n}Y_{j, n} \leq x s_n\right)- P\left (\sum_{j=1}^{m_n}Y_{j, n}^* \leq x s_n\right)\right|$$ 
$$+ C_{25}\;\int_{-T_2}^{T_2} \frac{1}{|t|}\left|E\left(e^{i\frac{t}{s_n}\sum_{j=1}^{m_n} Y_{j, n}^*}\right) - e^{-\frac{t^2}{2}}\right|dt+ C_{26}\;\frac{1}{T_2}$$
\begin {equation}
< \sup_{x \in R}\left|P\left (\sum_{j=1}^{m_n}Y_{j, n} \leq x s_n\right)- P\left (\sum_{j=1}^{m_n}Y_{j, n}^* \leq x s_n\right)\right|+ C_{25}(I_1 +I_2)  +  C_{26}\frac{1}{T_2}. \label{E5}
\end{equation}
where
$$ I_1 = \int_{-T_2}^{T_2}\frac{1}{|t|}\left|\prod_{j=1}^{m_n}\varphi_j(\frac{t}{s_n})-e^{-\frac{m_nt^2 s_{p_n}^2}{2s_n^2}}\right|dt$$ and 
$$I_2 = \int_{-T_2}^{T_2}\frac{1}{|t|}\left|e^{-\frac{m_nt^2 s_{p_n}^2}{2s_n^2}}-e^{-t^2/2}\right| dt.$$
The bounds for the  expressions on the right side of (\ref{E4}) and (\ref{E5}) are obtained from the Propositions \ref{P1},  \ref{P2}, \ref{P4},  Corollary \ref{C1} and the value of $T_2$.\\To obtain the final bound we compare $ n^{-q\alpha/8}, \; b_n^2\;n^{\theta(1-\alpha)-\alpha},\; b_n^{-1}n^{-\alpha/2}, \; n^{-3\alpha/8},$ $n^{-\alpha(q-2)/2}$ and $n^{(1-\alpha)\theta}$ for various values of $\alpha$ and $q.$ We consider the cases $2 < q \leq 3$ and $q > 3$ separately. \\ 
Let us say $ c_n > 0$ dominates over $d_n$ if $\frac{d_n}{c_n} \rightarrow 0$ as $n\rightarrow \infty$.	The  bound   contains terms some of which dominate over others. In the case $0 < q \leq 3$ the domination depends on the value of $q$ in the ranges $2 < q < 8/3$ and $8/3 \leq q \leq 3$ and the choice of   $\alpha$ in the ranges given below \\
$0 < \alpha <\frac{2\theta}{3+2\theta}$\\
$\frac{2\theta}{3+2\theta}< \alpha < \frac{2\theta}{q+2\theta}$\\
$\frac{2\theta}{q+2\theta} \leq \alpha \leq < \frac{\theta}{1+\theta}$\\
$\frac{\theta}{1+\theta} \leq \alpha < \frac{2\theta}{1+2\theta}$\\
$\frac{2\theta}{1+2\theta} \leq \alpha < 1. $\\\\ For $ 2 < q < 8/3$,
after a tedious but elementary analysis we get the bound

$$ C_{27} \frac{1}{ n^{ \alpha(q-2)/2}}I\left(\left(0 < \alpha \leq \frac{2\theta}{q+2\theta}\right) \cup \left(\frac{\theta}{1+\theta} < \alpha \leq \frac{2\theta}{q-2+2\theta}\right) \right) $$
$$+ C_{28}\frac{b_n^2}{n^{\theta - \alpha(1+\theta)}}	I\left(\frac{2\theta}{q+2\theta} \leq \alpha < \frac{\theta}{1+\theta}\right)
+C_{29} \frac{1}{n^{\theta(1-\alpha)}}I\left(\frac{2\theta}{q-2+2\theta} \leq \alpha < 1\right)$$
This can be simplified further. Since $n^{-\alpha(q-2)/2}$ decreases as $\alpha$ increases the best rate  contributed by the first term is for the maximum value of $\alpha$.  So we compare for $ \alpha =\frac{2\theta}{q+2\theta}$ and $ \frac{2\theta}{q-2+2\theta}$ and get the best rate  $ n^{-\frac{\theta(q-2)}{q+2\theta}}$. On the other hand for the same value of $q$, the second term gives the rate $b_n^2\;n^{-\frac{\theta(q-2)}{q+2\theta}}$, which is dominated by the previously obtained rate because $b_n^2 \rightarrow 0$ while the third term gives the rate $n^{-\frac{\theta(q-2)}{q-2+2\theta}} $ which too is dominated by  $n^{-\frac{\theta(q-2)}{q+2\theta}}$.  Thus for $2 < q \leq 8/3$ we get the rate  $n^{-\frac{\theta(q-2)}{q+2\theta}}$.\\
In the case $8/3 \leq q \leq 3$ the bound for the expression on the right side of (\ref{E4})  turns out to be 

$$ C_{28}\frac{b_n^2}{n^{\theta - \alpha(1+\theta)}} I\left( \frac{8\theta}{q+8+8\theta} \leq \alpha < \frac{\theta}{1+\theta}\right)
+C_{29} \frac{1}{n^{\theta(1-\alpha)}}I\left( \frac{8\theta}{q+8\theta} \leq \alpha < 1\right)$$
$$+C_{30} \frac{1}{n^{q\alpha/8}}I \left(\left(0 < \alpha < \frac{8\theta}{q+8+8\theta}\right)\cup \left(\frac{\theta}{1+\theta} \leq \alpha < \frac{8\theta}{q+8\theta}\right)\right).$$
Thus in the case $ 8/3 \leq q \leq 3$ the best rate is $ n^{-\frac{q\theta}{8+q+8\theta}}$. 
In the case $q \geq 3$ the  bound for the expression on the right side of (\ref{E4})  turns out to be $$C_{32}\frac{1}{n^{3\alpha/8}}\;I\left(\left(0 < \alpha < \frac{8\theta}{11+8\theta}\right)\cup \left(\frac{\theta}{1+\theta} \leq \alpha \leq \frac{8\theta}{3+8\theta}\right)\right) $$
$$+ C_{29} \frac{b_n^2}{n^{(1-\alpha)\theta-\alpha}}\;I\left(\frac{8\theta}{11+8\theta} \leq \alpha < \frac{\theta}{1+\theta}\right) + C_{30} \frac{1}{n^{(1-\alpha)\theta}}\left(\frac{8\theta}{3+8\theta} \leq \alpha < 1\right).$$
Thus in the case $q \geq 3$ the best rate is  $C_{33} \;n^{-\frac{3\theta}{11+8\theta}}$.\\

The best bound turns out to be  $$C_{22}\;n^{-\frac{3\theta}{11+8\theta}}$$ establishing the result.\\ This completes the proof of the Theorem \ref{T1}.
\begin{rem}
1. In  the border case of $q=8/3$ both the bounds $ n^{-\frac{\theta(q-2)}{q+2\theta}}$ and $n^{-\frac{q\theta}{q+8+8\theta}}$ coincide.\\
2. The bound in the Theorem \ref {T1} is independent of $\alpha$. However in \c{C}a\v{g}in \textit{et al.}  (2016) the bound depends on $\alpha$.\\
3.For $q=3$ the rate is $n^{-\frac{3\theta}{11+8\theta}}$, which as $ \theta \rightarrow \infty$, goes to $n^{-3/8}$ and this is far better rate than the rate $ n^{-1/5}$ given in Oliveira's book (2012). \\
4. As is to be expected the rate of convergence in the CLT improves as $q$ increases in the interval (2, 3). Further as  in the case of independent and identically distributed rvs the rate remains the same with finiteness of  the  moments of order $\geq 3$. \\ 
5. If $\mu_n$ is chosen as $e^{-\mu\alpha},\;\; 0 < \mu < 1/2$ instead of the above choice, the calculations become more complicated and we have to consider three cases; \textit{viz.,} $2 < q < \frac{1}{\mu}, \;  \frac{1}{\mu} < q < \frac{2\mu}{1-2\mu} $ and $\frac{2\mu}{1-2\mu} < q \leq 3 $ instead of $2 < q < 8/3$ and $8/3 \leq q < 3 $ when $q < 3$. The best rate turns out to be $ n^{-\frac{\mu\theta}{\mu+1+\theta}}$ for any choice of $q \in \left[\frac{2\mu}{1-2\mu},  3\right].$  Interestingly the above interval collapses to the single point set consisting of 3 when $\mu=3/8$.\\
  
\end{rem}
\section{Moderate deviation result} \c{C}a\v{g}in \textit{et al.}  (2016) recently obtained a moderate deviation result for associated rvs under strong conditions. Before we state and prove the moderate deviation result, we shall recall a result of Frolov (2005) and apply it to coupling block rvs introduced earlier. 

\begin{thm} (Theorem 1.1 in Frolov, 2005) \label {T2} Let $\{Y_{k, n}, k=1, 2, \ldots,k_n, n=1, 2, \ldots\} $ be an array of column-wise independent centered rvs with $EY_{k, n}^2=\sigma_{k, n}^2 < \infty$. Denote $T_n =\sum_{k=1}^{k_n}Y_{k, n}$ and $B_n=\sum_{k=1}^{k_n}\sigma_{k, n}^2$. Assume for some $q > 2, \;\;\;  E[Y_{k, n}^q\;I(Y_{k, n} > 0)]=\beta_{k, n} < \infty$, $B_n \rightarrow \infty$ and set $$M_n=\sum_{k=1}^{k_n} \beta_{k,n },\;\;\; L_n=\frac{M_n}{B_n^{q/2}},$$
$$\Lambda_n(t, s, \delta)=\frac{t}{B_n}\sum_{k=1}^{k_n}E(Y_{k, n}^2 I(-\infty < Y_{k, n} < -\delta \sqrt{B_n}/s)).$$ 	 Assume that $L_n \rightarrow 0$, and that for  each $\delta > 0,\; \Lambda_n (x^4, x^5, \delta) \rightarrow 0.$  If $x_n \rightarrow \infty$ such that
\begin{equation} \label {E6}
x_n^2-2 \log (1/L_n) -(q-1)\log \log (1/L_n) \rightarrow -\infty,
\end{equation}
then 
$$P(T_n \leq x_ns_n)=(1-\Phi(x_n)(1+o(1)).$$ 
\end{thm}	
Let the assumption $A_1$ hold for the original rvs $X_n$. Recall that the block rvs $Y_{k, n}^*,\; k=1, 2, \ldots, m_n$ are independent and identically distributed for each $n$ with $E|Y_{k,n}|^q < \infty$ where $q > 2.$ With $Y_{k, n}=Y_{k, n}^*, \;\;\; k_n=m_n,$ $$B_n = m_n\; s_{p_n}^2 \sim n \sigma^2 ,\;\;\; M_n \leq  m_n\; E|Y_{k, n}^*|^q  \sim n^{\alpha +(1-\alpha)q/2},\;\;\; L_n\sim n^{\alpha(2-q)/2}\rightarrow 0$$ as $n \rightarrow \infty$. Further 
\begin{equation}
\Lambda_n (x^4, x^5, \delta)  \leq  \frac{ x^4}{\sigma^2\;n^{1-\alpha}}E(Y_{k, n}^2 I(D_n)). \label {E7}
\end{equation}
where $D_n$ is the event $| Y_{k. n} |> \delta \sqrt{n}\sigma/x^5$ 
since $Y_{k, n}^*\stackrel{D}{=}Y_{k, n}.$ By the  H$\ddot{o}$lder's inequality and finiteness of moment of order $q$ for $Y_{k, n}$, we get 
$$E(Y_{k, n}^2 I(D_n)) \leq \left( E(Y_{k, n}^2 I(D_n))^{q/2}\right)^{q/2}\left(E(I(D_n))^{q/(q-2)}\right)^{(q-2)/q}$$
$$\leq (E|Y_{1, n}|^q)^{2/q}(P(D_n))^{(q-2)/2} \leq p_n \left( \frac{E|Y_{1, n}|^q x^{5q}}{\delta^q n^{q/2}\sigma^q}\right)^{(q-2)/q},$$
which results in the following bound from (\ref{E7})
$$\Lambda_n (x^4, x^5, \delta) \leq x^4 \left(\frac{p_n^{q/2} x^{5q}}{\delta^q\sigma^q n^{q/2}}\right)^{(q-2)/q}$$
$$ \leq C_{34}\frac{x^{5q-6}}{n^{\alpha(q-2)/2}}.$$
If $x=x_n \sim  (\log n)^\kappa$ and  $ \kappa > 0$  we then have $\Lambda_n (x_n^4, x_n^5, \delta) \rightarrow 0$ as $n \rightarrow \infty, $ so that all the conditions of the Theorem \ref{T2} hold and we then get the following moderate deviation result for the coupling block rvs ${Y_{k, n}^*}.$ 	

\begin{thm} \label {T3}
If  $\{X_n\}$ 	is a sequence of centered associate rvs satisfying the assumption $A_1$ then for the coupling block rvs $Y_{j, n}^*$ 
$$ P\left(\sum_{j=1}^{m_n}Y_{j, n}^* > x_n s_n\right)=(1-\Phi(x_n))(1+o(1))  $$
whevever $x_n$ satisfies $$\limsup_{n\rightarrow \infty}\frac{x_n^2}{\log n}= \lambda < \alpha (q-2).$$ 
\end{thm}
\begin{rem} In the Theorem 4.2 of \c{C}a\v{g}in \textit{et al.}  (2016) the Assumption (B2) states the condition differently but a close look at the proof reveals that they indeed use  $\limsup_{n\rightarrow \infty}\frac{x_n^2}{\log n} < 1$ which is similar to our assumption.
	\end {rem}
\begin {cor} \label {C2}
Recall $\mu_n=n^{-3\alpha/8}.$ If $x_n$ satisfies the relation (\ref{E6}) then so will $x_n \pm \mu_n$ and we have 
$$ P\left(\sum_{j=1}^{m_n}Y_{j, n}^* > (x_n \pm \mu_n) s_n\right)=(1-\Phi(x_n))(1+o(1))  $$
because $\mu_n =o(1-\Phi(x_n)).$ Here we use the fact $|\Phi(x+\epsilon) - \Phi(x)| <\epsilon$.
\end{cor}
Now we state and prove the moderate deviation result for $S_n.$
\begin{thm} \label {T4}
Let $\{X_n\}$ be a sequence of centered stationary associated rvs satisfying the  assumptions $A_1$ and $A_2$. Assume further \\
(i) $ lim\sup_{n\rightarrow \infty}\frac{x_n^2}{\log n} = \lambda < \frac{q-2}{2}$,\\
(ii) $ \theta$ in the assumption $A_2$ is such that
\begin{equation}\label{E8}
\theta > 1+\lambda.
\end{equation} Then 
\begin{equation} \label {E9}
P(S_n > x_n\; s_n)=(1-\Phi(x_n))(1+o(1)).  
\end{equation}
\end{thm}
\textbf{Proof}
Choose $\alpha$ in the definition of $p_n$ such that 
\begin{equation}\label {E10}
\frac{1}{2} < \alpha < \frac{2\theta-\lambda}{2\theta+2}.
\end{equation}
This is possible because of the assumption at  (\ref{E8}). Let $ \epsilon_n= n^{-\epsilon}$ where 
\begin{equation}\label{E11}
0 < \epsilon < \frac{q\alpha-\lambda}{2\;q}.
\end{equation}
This is possible because $ \lambda < (q-2)/2 $ and $ \alpha > 1/2.$ The stated result follows from the Corollary \ref{C2} and the assumption (i) above if we prove
$$ (a)\;\;\;\left|P(S_n > x_n\;s_n)-\left(\sum_{j=1}^{m_n}Y_{j, n} > (x_n\pm \epsilon_n)\;s_n\right)\right| = o(1-\Phi(x_n)) $$
and 
$$(b)\;\;\;\left|\left(\sum_{j=1}^{m_n}Y_{j, n} > (x_n\pm \epsilon_n)\;s_n\right) - \left(\sum_{j=1}^{m_n}Y_{j, n}^* > (x_n\pm \epsilon_n)\;s_n\right)\right|= o(1-\Phi(x_n)).
$$
To prove (a) recall from the Proposition \ref{P1}
$$ \left|P(S_n > x_n\;s_n)-P\left(\sum_{j=1}^{m_n}Y_{j, n} > (x_n\pm \epsilon_n)\;s_n\right)\right| \leq P\left(|Y_{m_n+1, n}| >\epsilon_n\;s_n \right)$$
\begin{equation}\label{E12}
< C_{35}\; \frac{p_n^{q/2}}{\epsilon_n^q\; n^{q/2}} < C_{36}\; n^{-q\;(\alpha-2\epsilon)/2}
\end{equation}
We get the result (a) if $$\frac{\sqrt{\log n}}{n^{(q(\alpha-2\epsilon)- \lambda)/2}}\rightarrow 0$$ which follows from (\ref{E10}).
\\Next to prove (b) recall from  the Proposition \ref{P2}
$$\left|\left(\sum_{j=1}^{m_n}Y_{j, n} > (x_n\pm \epsilon_n)\;s_n\right) - \left(\sum_{j=1}^{m_n}Y_{j, n}^* > (x_n\pm \epsilon_n)\;s_n\right)\right| $$
\begin{equation}\label{E13}
< C_3\;\frac{b_n^2}{n^{\theta-\alpha(1+\theta)}}I\left(\frac{2\theta}{3+2\theta} \leq \alpha < \frac{\theta}{1+\theta}\right)+C_4\;\frac{1}{b_n n^{\alpha/2}}I\left(\alpha < \frac{2\theta}{3+2\theta}\right).
\end{equation}
The first term on the right side above is $o(1-\Phi(x_n)) $ because (\ref{E10}) implies $\theta -\alpha (1+\theta) > \lambda/2$. The second term on the right side of (\ref{E13}) is $o(1-\Phi(x_n))$ because  $ \lambda < \alpha.$ This completes the proof of the Theorem.
\begin{rem} 
\c{C}a\v{g}in \textit{et al.}  (2016) proved the Theorem \ref{T2} making complicated assumptions of the type $A_2$ as well as $A_3$ with the conditions that  $\theta > 4 $ and $ q > 3.$  Further our proof does not require dealing with odd numbered and even numbered blocks separately nor does it need introduction of Gaussian centered varibles similar to odd and even block sums.\\
\textbf{Acknowledgement} The author thanks Professor B. L. S. Prakasa Rao for bringing the paper by  \c{C}a\v{g}in \textit{et al.}  (2016) to his notice.
\end{rem}

\textbf{ References} \\

\noindent \bf{Birkel, T.} (1988) On the convergence rate in the central limit theorem for associated processes, \textit{Ann. Probab}, 16, 1689-1698.\\
\bf{\c{C}a\v{g}in, T., Oliveira, P. E. and Torrado, N.} (2016) A moderate deviation for associated randon variables, \textit{Jl. Korean Statist. Soc.}, 45, 285-294.\\
\bf{Dembo,A.  and Zeitouni, A.} (1998) \textit{Large deviation techniques and applications}, 2nd Edn., Springer, New York.\\
\bf{Frolov, A. N. }(2005) On probabilities of moderate deviations of sums for independent random variables, \textit{Jl. Math. Sci.}, 127, 1787-1796.\\
\bf{Hollander, F. D.} (2000) \textit{Large deviations}, Fields Inst. Monographs, Amer. Math. Soc., Providence, Rhodes Island.\\
\bf{Newman, C.M.} (1980) Normal fluctuations and FKG-inequality, \textit{Commun. Math. Phys.}, 74, 119-128. \\
\bf{Oliveira, P. E.} (2012)\textit{ Asymptotics for associated random variables}, Heidelberg, Springer.\\
\bf{Tikhomirov, A. N.} (1980) On the convergence rate in the central limit theorem for weakly dependent random variables, \textit{Theory Probab. Appl.}, 25, 790-809.\\
\bf{Vardhan, S. R. S.} (1984) \textit{Large deviations and applications}, S. I. A. M., Philadelphia.\\
\bf{Wang, W.} (2015) Large deviations for sums of random vectors attracted to operator semi-stable laws,\textit{ Jl. Theor. Probab. DOI 10.1007/s10959-015-0645-5.}\\
\bf{Wood, T. E.} (1983) A Berry-Esseen theorem for asociated random variables, \textit{Ann. Probab.,}
11, 1041-1047.\\\\

\end{document}